\documentclass{article}
\usepackage{graphicx}
\usepackage{amsmath} \usepackage{euscript} \usepackage{amssymb}
\usepackage{amsthm} \usepackage{amsopn} 

\begin{document}    
     
\newtheorem{prop}{Proposition}
\newtheorem{defi}[prop]{Definition}
\newtheorem{thm}[prop]{Theorem}
\newtheorem{ex}[prop]{Example}
\newtheorem{lem}[prop]{Lemma}
\newtheorem{rem}[prop]{Remark}
\newtheorem{cor}[prop]{Corollary}
\newtheorem{conj}[prop]{Conjecture}
\title{ \vspace*{-2em}$S_k$-Holonomy on Coloring Complexes of $M^n$ \\
                                            with  Applications to the \\
                                Poincar\'{e} Conjecture and $4$-Color Theorem}
\author{Daniel Kling\thanks{email: dkling@foldedstructures.com}}
\date{May 22, 2017}
\maketitle
\large
\begin{abstract} 
  A natural class of coloring complexes $X$ on closed manifold $M^n$
is investigated that gives a holonomy map $\mbox{Hol}_X: \pi_1(M) \to S_{n+1}$.
By a $k$-multilayer complex construction the holonomy map may be defined
to any finite permutation group
 $\mbox{Hol}_X: \pi_1(M) \to S_{n+k}$, $k>0$.
Under isotopy of $X$ and surgery on $B^n \subset M^n$
a holonomy class of complexes $[X]$ is defined with
$[X]=[Y] \iff  \mbox{Hol}_X = \mbox{Hol}_Y$.
It is also shown that for any homeomorphism $f:\pi_1(M) \to S_{n+1}$
there is a complex $X$ on $M$ with  $\mbox{Hol}_X =f$. 
These results are applied to express the $4$-color Theorem and the Poincar\'{e} Conjecture
as the existence and uniqueness, respectively, of a certain holonomy class. 
Several other applications are suggested.
\end{abstract}

\tableofcontents

\section{Introduction}    
 Throughout the paper $M^n$ will be a closed connected triangulizable and smooth manifold, 
with no boundary unless explicitly stated.
A complex $X$ on $M$ that is dual to a simplicial triangulation will be called a regular complex,
and written $X$ cpx on $M$.
These have regular degree in the sense that
for $1\leq k\leq n$ there are $(n-k+2)$
 $k$-cells meeting at each $(k-1)$-cell.
For example, a rack of soap bubbles in $2$ or $3$ dimensions has this property.
The top dimensional components of $M^n$ will be called regions.
For clarity the $i$-skeleton  of $X$  will be denoted $X^{(i)}$
and its set of $i$-cells by $X^{\langle i \rangle}$.

We start by considering $(n+1)$-colorings of the regions of $X$ cpx on $M^n$.
The regular degree of $X$ implies there is at most one such forced coloring.  
Conditions are shown for $X$ to be locally $(n+1)$-colorable,
and this yields a holonomy homomorphism $\mbox{Hol}_X: \pi_1(M)\to S_{n+1}$.
A holonomy class $[X]$ of cpx on $M$ is developed that is natural 
with respect to both surgery on $B^n$ and the holonomy homomorphism.
Moreover a $1-1$ correspondence is shown between the holonomy classes on $M^n$
and the homomorphisms $f:\pi_1(M) \to S_{n+1}$.

The method is extended to sets of $k$ cpx on $M$, where the layers are transverse.
Conditions for this to be locally $(n+k)$-colorable are given,
 and for these layered cpx there is a holonomy homomorphism $\pi_1(M) \to S_{n+k}$.  
With this vocabulary the   $4$-Color Theorem and the Poincar\'{e} Conjecture are 
compared by 
expressing them as special cases of similar statements.
Further work and applications for studying $\pi_1(M)$ and colorings on $M^n$ are mentioned.

\section{$S_{n+1}$ Coloring Holonomy on $M^n$} 
\begin{defi}\label{defi_col}
Let $X$ be a cpx on $M^n$.  For $c \in \mathbb{N} $, 
let $\{1,\ldots, c\}$ represent $c$ distinct colors, and let $f:X^{\langle n \rangle }\to \{1, 2,\ldots, c\}$
denote an assignment of colors to the set of regions of $X$.  $f$ is called a coloring
if for distinct regions $r_1, r_2$, we have    $ r_1 \cap r_2 \neq \emptyset \Rightarrow  f(r_1) \neq f(r_2)$.
Then say $f$ is a $c$-coloring of $X$. 
 If $\forall r \in X^{\langle n \rangle}$, $\exists$ an $\epsilon$-tubular neighborhood $N_\epsilon(r)$
such that
the neighborhood's  complex
$X|_{N_\epsilon(r)}$ is $c$-colorable then say $X$ is locally $c$-colorable.
\end{defi}

\begin{prop} \label{exists_cpx} Let $M^n$ be given.
\begin{enumerate}
\item  $\forall \epsilon >0, \exists (n+1)$-colorable cpx $X$ on $M$ with 
each cell of $X$ contained in a ball or radius $\epsilon$.
\item $\exists (n+1)$-colorable cpx $X$ on $M$ with $X$ having $(n+1)$-regions.
\end{enumerate}
\end{prop}

\begin{proof} (1) First choose a triangulation $T$ on $M^n$ with each cell of diameter  $< \epsilon/2$.
We will construct a cpx X with one region for each cell of $T$.
For the vertices of $T^{(0)}$ we may choose $\delta > 0$ with 
tubular neighborhood $U_0 = N_\delta(T^{(0)})$ to form a collection of non-intersecting balls.  
Next choose a smaller $\delta >0$ with tubular neighborhood $N_\delta(T ^{(1)})$ and put
$U_1$ equal to the connected components of
 $N_\delta(T ^{(1)}) \setminus \mbox{ Int}( U_0)$.
This may be chosen to be a collection of non-intersecting cylinders with ends meeting $U_0$.
For each $i < n$ we may  choose a yet smaller $\delta$ so that $U_i$
is a disjoint collection of tubular neighborhoods around the $i$-cells in
$X^{(i)} \setminus \mbox{Int}        \bigcup \{ U_0, \ldots ,  U_{i-1}  \} $ 
and put 
$U_n$ equal to the connected components of $ M \setminus \mbox{Int} \bigcup \{ U_0, \ldots , U_{(n-1)}    \}$.
Then 
  $\bigcup \{ U_{0}, \ldots ,  U_{n}   \} $
forms a collection of closed $n$-cells whose CW-complex satisfies the regular degree condition and is a cpx on $M$. 
 Coloring each region of $U_i^{\langle n \rangle}$ in color $i$
gives the claim.

For (2), we may choose $(n+1)$-colorable cpx $X$ on $M$ by (1).
Note each edge $e$ of $X$ is surrounded around its interior by $n$ mutually adjacent regions, 
with the regions at both ends in the remaining color.  
The apparently two end-regions of $e$  could be connected 
away from the edge and actually the same region.
Induction will be on the number of regions in $X$.
Let $c$ be a color.  Observe the union of the $c$-colored regions and 
the edges only adjacent to these regions on their ends forms a connected set.
Thus if there are multiple $c$-colored regions than one of these edges
 must join two distinct $c$-colored regions $r_1$ and $r_2$.  
Connect  $r_1$ and $r_2$ by tunneling along the edge.  Specifically, 
replace the edge with a tubular neighborhood $N$ and form a new region
as the union of  $r_1, r_2$ and $N$.  
The new region may be $c$-colored and yields a $(n+1)$-colorable cpx with one fewer regions.
Proceeding by induction we may reduce to a $(n+1)$-colorable cpx with one region of each color.
\end{proof}

\begin{rem}\label{convex}
It appears by choosing triangulation $T$ in Proposition~\ref{exists_cpx} to have geodesic simplicies, 
the  regions of $X$ may be constructed to have diameter less than $\epsilon$ and be convex.
\end{rem}

Next we define holonomy for locally $(n+1)$-colorable cpx on $M^n$.
Suppose $M^n$ has a cpx $X$, and let $W$ be a tubular neighborhood of 
$X^{(1)}$.
Let $v_0 \in X^{\langle 0 \rangle}$.  The degree of each vertex is $n+1$, 
so the regions in a neighborhood of $v_0$ may be $(n+1)$-colored.
Let $\gamma = (v_0, \ldots, v_k)$ be a path on $X^{(1)}$. 
The coloring around $v_0$ may be extended to the tubular neighborhood of 
$\overline{v_0,v_1}$ in $W$.  
Since the degree of each edge is $n$, one color is not used along 
Int$(\overline{v_0,v_1})$.
At $v_1$ however, there is another adjacent region in $W$, and this is forced to be assigned the unused color.
Likewise the $(n+1)$-coloring around $v_1$ extends along 
$\overline{v_1,v_2}$ 
in $n$-colors, and forces a unique $(n+1)$-coloring around $v_2$.
If $v_0 = v_k$ and so $\gamma$ is a cycle, then for initial choice of 
$(n+1)$ colors at $v_0$, the coloring at $v_k = v_0$ 
gives another coloring and a permutation in $S_{n+1}$.
The composition of cycles in $\pi_1$ corresponds to the composition of permutations.
This defines the holonomy map from $\pi_1(X^{(1)}) \to S_{n+1}$.  This extends uniquely to $\pi_1(W)$.

\begin{defi}\label{defi_HolW} For $X$ a cpx on $M^n$, and  $W$ a tubular neighborhood of $X^{(1)}$, 
 let the holonomy map be defined as above and be denoted by
\[ \mbox{Hol}_W: \pi_1(W) \to S_{n+1}  \]
\end{defi}

If $b \in  X^{\langle 2 \rangle}$ is an even sided $2$-cell, then the cycle 
$\gamma = \mbox{bdry}(b)$ yields $\mbox{Hol}_X(\gamma) = 1$.  
To see this note the degree of each edge of $b$ is $n$, with $n-1$ of the regions adjacent to Int$(b)$.  
Thus there are only two colors available for regions skirting around $b$,
that is those regions adjacent to $\gamma$ and not adjacent to the interior of $b$.
These perimeter regions must then alternate in color, and $b$ even-sided yields the identity holonomy.
It also shows when a $2$-cell $b$ is even-sided,  a neighborhood of $b$ is $(n+1)$-colorable.
Thus

\begin{lem}\label{even_sided} Let $X$ be a cpx on $M^n$.  Then the following are equivalent
\begin{enumerate}
\item   All $2$-cells of $X^{\langle 2 \rangle}$ are even sided
\item   $\mbox{Hol}_{X^{(1)}}: \pi_1(X^{(1)}) \to S_{n+1}$ extends to 
 $\mbox{Hol}_{X^{(2)}}: \pi_1(X^{(2)}) \to S_{n+1}$
\item  An $\epsilon$-nbd of all $2$-cells of $X^{(2)}$ is $(n+1)$-colorable.
\end{enumerate}
\end{lem}

\begin{prop}\label{HolX}
For $X$ cpx on $M^n$ 
with all $2$-cells even sided, then
the holonomy extends uniquely from
 $\mbox{Hol}_{X^{(1)}}: \pi_1(X^{(1)}) \to S_{n+1}$ to
\[ \mbox{Hol}_X: \pi_1(M) \to S_{n+1}\]
\end{prop} 
\begin{proof}
Let $X$ satisfy the hypothesis and cycle $\gamma$ be given with base point $v_0$ in  $X^{\langle 0 \rangle}$. 
Then $\gamma$ may be homotopied to $\gamma' \subset  X^{(1)}$. 
If $\gamma$ is also homotopied to $\gamma'' \subset  X^{(1)}$, 
then one may show  $\gamma' $ and  $\gamma''$ are homotopic in  $X^{(2)}$.  But all homotopies in $X^{(2)}$
are compositions of homotopies across even-sided $2$-cells, so $\gamma' $ and  $\gamma''$
have the same holonomy.  Thus
$\mbox{Hol}_X: \pi_1(M)  \to S_{n+1}$
is uniquely defined.
\end{proof}

\begin{thm}\label{loc123}
Let $X$ be a cpx on $M^n$.  Then
\begin{enumerate}
\item        $X$ is locally $(n+1)$-colorable $\iff$ All $2$-cells in  $X^{(2)}$ are even sided.
\item If $M$ is orientable then:
               $X$ is locally $(n+1)$-colorable $\iff$  $X^{(1)}$ is even cyclic
\item If $M$ is simply connected then:
               $X$ is locally $(n+1)$-colorable $\iff$  $X$ is $(n+1)$-colorable.
\end{enumerate}
\end{thm}
\begin{proof}
(1)  was shown in Lemma~\ref{even_sided}.  
For (2),  $X^{(1)}$ even cyclic $\Rightarrow$ all $2$-cells are even-sided
$\Rightarrow$ $X$ is locally $(n+1)$-colorable.
For the converse, $X$ orientable implies there is a universal orientation from the tangent field at the base point.
Thus we may determine a universal orientation for each $n+1$ color assignment around the vertices,
 based on the sign of the determinant of the
color-matching approximating linear transformation at the vertex.
Since this alternates at the ends of each edge, all cycles must be even length.
For (3), clearly $(n+1)$ colorable $\Rightarrow$ locally $(n+1)$ colorable.  
If $M$ is simply connected, 
let base point $v_0$ and vertex $v_1$ with paths $\gamma, \gamma'$ from $v_0$ to $v_1$ be given. 
By holonomy both paths give a coloring at $v_1$.
Since 
$\mbox{Hol}_X( \gamma\gamma'^{-1}) =\mbox{Hol}_X( 1) =1 \in S_{n+1} $,
$\gamma$ and $ \gamma'$ have the same holonomy, and
the coloring at each vertex 
can be constructed independently from the choice of path from the base point, 
and $X$ is $(n+1)$-colorable.
\end{proof}

\begin{defi}
Let $X$ be a cpx on $M^n$.  Say $X$ is an ecpx on $M$ if all $2$-cells of $X$ are even sided.
\end{defi}
So $X$ ecpx on $M^n \iff  X \mbox{ locally } n+1 \mbox{ colorable }
\iff  
\mbox{Hol} \mbox{ extends from } X^{(1)} \mbox{ to }X$.

\section{Fully Transverse Submanifolds}                         

\begin{defi}\label{bdry_cpx}
Let $M^n$ with boundary $N^{n-1}$ be given and suppose a complex $X$ is given on $(M,N)$
with embedded cells, $X$ regular degree in the interior of $M$ and $X|_N$ a regular degree complex on $N$.  Then say 
$X$ is a cpx on $(M,N)$.
\end{defi}
For example for $M=B^3$, the vertices in the interior have edge-degree $4$ while those on $S^2$ have edge-degree $3$. 
So cpx with boundary are dual to simplicial triangulations with boundary.
Similarly for hyper-submanifolds we have
\begin{defi}\label{cpx_section}
Let $M^n$ be closed with possible boundary, and  $X$ a cpx on $M$. 
 Let $W^{n-1}$ be a submanifold that meets $X$ transversely.
Suppose further, $\forall i \leq n$, and $\forall b \in X^{\langle i \rangle }$ with 
$W \cap b \neq \emptyset$,
$W \cap b \cong D^{i-1}$.
%
%
Then say $W$ is a fully transverse hypersurface of $X$.
\end{defi}

Note when $X$ is a cpx on $(M,N)$, for $X^{(1)}$ to be considered even cyclic or for 
 $X^{(2)}$ to have even-sided $2$-cells, only the cells of  $X$ missing $N$ are required to satisfy the condition.
If $f:X^{\langle n \rangle} \to \{1, \dots , n+1\} $ is a coloring, then this gives a coloring 
$f|_W : W^{\langle n-1 \rangle} \to \{1, \dots , n+1\} $.
Also, a coloring of $X^{\langle n \rangle}$ gives a coloring of $(X\cap N)^{\langle n-1 \rangle}$
where the colors of the boundary hyper-regions inherit the color of their unique adjacent region of $X$.
If  $f : W^{\langle n-1 \rangle} \to \{1, \dots , n+1\} $ is a coloring, a coloring $F:X^{\langle n \rangle} \to \{1, \dots , n+1\} $ s.t. $F|_W = f$
is called an extension of $f$.  

\begin{prop}\label{fill_from_bdry}
Let $(M^{n},N^{n-1})$ be given with  $Y$ an $(n+1)$-colorable cpx on $N$.
Then $Y$ has an extension to an $(n+1)$-colorable ecpx $X$ on $M$.
\end{prop}
\begin{proof}
Assume any $k$-coloring of $S^n$, $k \geq n+2$, can be extended to a $k$-coloring of $B^n$.  
A tubular neighborhood of $N$ may be taken providing a partial complex $Y \times [0,\epsilon] \subset M$.
The remainder of $M$ may be subdivided 
into a relatively coarse complex $Z$ with cells meeting $Y \times \epsilon$  fully transversely.
We may introduce relatively small $n$-cells 
centered at any vertices of $Z$ that lie in
$M \backslash (Y \times [0,\epsilon])$.
By hypothesis these and the $(n+1)$-colored cells of $Y \times [0,\epsilon]$ may be extended to a 
relatively fine cpx covering each of the $1$-cells of $Z$.
This gives a partial coloring of $M$ that covers the boundary of each of the $2$-cells of $Z$.
By hypothesis the partial coloring may be extended with relatively fine $n$-cells to cover the $2$-cells of $Z$.
By induction the $(n-1)$-skeleton of $Z$ is covered with colored cells and lastly the remaining portion of the
$n$-cells of $Z$ can be filled in with colored $n$-cells, producing a full coloring of $M$ that extends the coloring on $N$.  
The proof is complete with the following lemma.
\end{proof}

\begin{lem}\label{Sn_fill}
Let $W$ be a cpx on $S^n$ with coloring $f:W\to\{1, \ldots, n+2\}$.  Then there exists a cpx $X$ on $B^{n+1}$
with $X|_{S^n} = W$ and coloring $F:X \to \{1, \ldots,n+2 \}$ that extends $f$.
\end{lem}
\begin{proof}
Clear for $n=0$.  For $n=1$ let a partition $W$ of $S^1$ and coloring of $W$ be given.  
Consider $S^1$ to lie in the plane, and thicken $W$ radially to form an annulus of colored $2$-cells.
If the color regions use only $2$ colors, use the third color in the center to give the $3$-coloring of $B^2$.
If the inside edge of the annulus has $3$ colors but one color is used only once,
expand this region to fill the annulus and give a $3$-coloring of $B^2$.
Otherwise, some region $r_1$ on the inside of the annulus is bounded by two regions $r_2$ and $r_3$ of distinct colors.
Then $r_2$ and $r_3$ may be grown to meet and cover $r_1$ so that it is no longer exposed on the inside edge of the annulus.
Proceeding by induction reduces the number of exposed interior faces and fills the disk.

For $W$ a cpx on $S^2$ the argument is similar.  $W$ is thickened to make $4$-colored 3-cells.  
These are grow to fill $B^3$ as follows.  If there is only one blue region, 
it may be grown to fill the inside of the shell and give $X$ on $B^3$ with a  $4$-coloring $X$.
Otherwise an interior blue face may be selected.
The boundary of the face is $S^1$ and surrounded by at most $3$-colors.  
Applying the $n=1$ case extends the boundary regions to meet and cover over the blue face.
By induction we may reduce to having only one exposed blue face, 
which is then grown to fill the interior and provide a $4$-coloring of $B^3$.

By induction on $n$ the argument applies to $S^n$.
\end{proof}

This gives the following reformulation and extension of the $4$-Color Theorem\cite{4C}.

\begin{cor}\label{Sn_even}  
Let W be a cpx on $S^n$.  $W$ is $(n+2)$-colorable
$ \iff \exists$ an ecpx $X$ on $B^{n+1}$ with  $X|_{S^n} = W $.
\end{cor}
\begin{proof}
Since $B^{n+1}$ is simply connected, $X$ ecpx
 $\Rightarrow X$ is $(n+2)$-colorable 
$\Rightarrow W=X|_{S^n}$ is $(n+2)$-colorable.
Conversely, 
$W \ (n+2) \mbox{-colorable } \Rightarrow \exists (n+2) \mbox{-colorable}$ extension  $X $ on $ B^{n+1}  \Rightarrow X$ 
is even cyclic.
\end{proof}

\section{Holonomy Class}          

For ecpx $X$ on manifold $M^n$ it is helpful to define a local surgery 
that does not impact the holonomy of $X$. 
Let $B^n \subset M^n$ be given with boundary $S^{n-1}$ fully transverse to $X$.
We wish to replace ecpx $Y=B^n \cap X$ with a choice of $Y'$ ecpx on $B^n$,
where it is required
$S^{n-1}$ is fully transverse to $Y'$,  
$S^{n-1}  \cap Y' = S^{n-1} \cap Y$, and
$Y' \cup_{S^{n-1}} Y$ is an ecpx on 
$B^n  \cup_{S^{n-1}} B^n = S^{n}$.
Then
keeping $S^{n-1}$ fixed we may cut out $Y$ and sew in $Y'$.
Note each partial cell cut out of $Y$  is re-spliced with a partial cell of $Y'$.
Call this a holonomy surgery on $B^n$.

\begin{defi}\label{Hol_Xi}
Let $(M^n, N^{n-1})$ be a manifold with possible boundary  $N$.  Let $X$ be ecpx and $X'$  cpx on $(M,N)$
that agree on an $\epsilon$ neighborhood of $N$.
If
there exists a finite number $k$ of cpx's $X_i$ on $M$ with $X_0 = X$ and $X_k =X'$,
and
with $X_i \cap N_\epsilon = X \cap N_\epsilon$
 s.t.
$\forall i \leq k$,
  either 
\begin{enumerate}
\item $X_i$ and $X_{i+1}$ are ambient isotopic keeping $N_\epsilon$ fixed.
\item $X_{i+1}$ is obtained from a holonomy surgery  on $X_i$  keeping $N_\epsilon$ fixed.
\end{enumerate}
Then write 
$X \approx X'$, and define $[X]$ to be the collection of all cpx $X'$ s.t. $X\approx X'$.
\end{defi}

\begin{prop}\label{approx_to_Hol} Let $X$ be ecpx on $M^n$ and $Y$ a cpx on $M$.
 Then $X \approx Y$ implies  $Y$ is an ecpx on $M$ and moreover
$\mbox{Hol}_X =  \mbox{Hol}_Y$.
\end{prop}

\begin{proof}
Suffices to prove the case $k=1$ of Definition~\ref{Hol_Xi}.  
If $X$ and $X'$ are ambient isotopic relative to $L_\epsilon$, the claim holds.
For the case where $X$ and $X'$ differ by a surgery on $B^n$
put $Y = B \cap X$ and $Y' = B \cap X'$ and let $S^{n-1}$ denote the bdry$(B)$.
We need to confirm all $2$-cells are even sided after surgery.
Let $r' \in X'^{\langle 2 \rangle}$ be given.  
Clearly if $r' \cap S^{n-1} = \emptyset$ then $r'$ is even-sided.
If $r' \cap S^{n-1} \neq \emptyset$, let $r$ be the $2$-cell of $X$
with $r \setminus B = r' \setminus B$.
So $r^{\langle 0 \rangle } \setminus B = r'^{\langle 0 \rangle } \setminus B$, and 
 $r^{\langle 0 \rangle } \cap B + r'^{\langle 0 \rangle } \cap B$ is even.
The latter is because the surgery used an ecpx on $S^n$.
Thus the parity of   $r^{\langle 0 \rangle }$ equals the parity of $ r'^{\langle 0 \rangle }$,
and $r'$ is even sided.

To see  $ \mbox{Hol}_X(\gamma)= \mbox{Hol}_{X'}(\gamma),  \forall \gamma \in \pi_1(M)$,
one can either homotopy $\gamma$ to miss $B$ or observe two colorings on the regions intersecting $S^{n-1}$
agree on $X$ and $X'$ so paths entering and emerging from $B$ would have the same holonomy in $X$ and $X'$.

\end{proof}

To prove the existence in Proposition~\ref{exists_cpx} extends to holonomy equivalence  we will need an interesting lemma.

\begin{thm}\label{small_cells}
Let $X$ ecpx on $(M^n, N^{n-1}),  \epsilon >0$.  Then 
\begin{enumerate}
\item
$\exists X$ ecpx on $(M,N)$ with $X' \approx X$
and all cells of $X$ having diameter $< \epsilon$.
\item
$\exists X$ ecpx on $(M,N)$ with $X' \approx X$
and $X'$ having $(n+1)$ regions.
\end{enumerate}
\end{thm}

\begin{lem}\label{PQ}
The following logical statements $P(n)$ and $Q(n)$ are true $\forall n > 0$: \\
\vspace{-2pt}

\noindent 
$P(n) \iff $ \parbox{.86\textwidth}{Let $X$ be a $(n+2)$-colored cpx on $S^n$, $\epsilon > 0$.  Then $\exists$ an $X'$ cpx on $S^n$,
                            with $X'$ fully transverse to $X$, $X'$ $(n+2)$-colorable, no like color regions of $X$ and $X'$ intersect,
                            and the cells of $X'$ have diameter $< \epsilon$.}
\vspace{10pt}

\noindent
$Q(n) \iff $ \parbox{.86\textwidth}{Let $X$ be a $(n+2)$-colored cpx on $S^n$, $\epsilon > 0$.  Then $\exists$  
                            $(n+2)$-colored $X'$ ecpx on $B^{n+1}$,
                            that is fully transverse to $X$, with 
                            intersecting regions of $X^{\langle n \rangle}$ and $X'^{\langle n+1 \rangle}$ having distinct colors,
                            and all cells of $X'$ having diameter $< \epsilon$.}
\end{lem}
\begin{proof} (Lemma)
 The condition $Q(0)$ is also true and asserts that any $2$-coloring of $S^0$ extends to an arbitrarily fine $2$-coloring of $B^1$.
To show $P(1)$, let $3$-coloring $X$ of $S^1$ be given, and $\epsilon > 0$.
Construct $X'$ by first covering the vertices with non-intersecting $B^1$ of radius $< \epsilon$.  
These may be $3$ colored distinct from $X$.  
Each edge of $X$ has its ends covered by at most two of the non-edge colors,
 so by $Q(0)$ the $2$-colored covering may be completed across the edge with all cells having diameter $ < \epsilon$.

We will next show $P(n) \Rightarrow Q(n)$.  Let $X$ be a $(n+2)$-coloring of $S^n$,
 and assume  $S^n$ is the unit sphere in $\mathbb{R}^{n+1}$.
 Let $\epsilon >0$.  Choose $(n+2)$-coloring $X'$ fully transverse to $X$, with all cell diameters $< \epsilon / 2$.
We may form $X' \times [0, \epsilon /2]$ by extending $X'$ radially inward by $ \epsilon / 2$. 
Each of the regions has diameter $< \epsilon$.
The inner boundary of this is $(n+2)$-colored sphere of radius $1 - \epsilon$.
Repeating the procedure with successive shells
and applying Lemma~\ref{Sn_fill} to the last shell
we may fill $B^{n+1}$ with $(n+2)$-colored regions of diameter $ < \epsilon$.
For $S^n \subset B^{n+1}$ not Euclidean, we may choose smooth map from the Euclidean metric to $S^n$ with Lipschitz number $h$.
By picking the $X'$ to have cells of diameter $< \epsilon / (2h)$ the argument follows.

Next let $n >0$.  Suppose $\forall m < n, Q(m)$.  If it is shown $P(n)$, then the proof is complete.
Let $X$ be $(n+2)$-colored on $S^n$.  Let $\epsilon >0$.  We may choose disjoint  $B^n \subset N_\epsilon(x)$ for 
$x \in X^{(0)}$.  There is $1$ choice of the $n+2$ colors for each $B$.  
The interior of each edge in $X^{\langle 1 \rangle}$ is adjacent to all but $2$ of the colors in $X$,
so by $Q(0)$ it may be $2$-colored with segments of length $< \epsilon /2$.  A small tubular neighborhood
of the segment gives regions of diameter $<\epsilon$.  Proceeding by induction we may cover all of the edges.

Suppose we have covered all of $X^{(i)}$ with regions of diameter $ < \epsilon$. 
Let $b \in X^{\langle i+1 \rangle}$.  The interior of $b$ is adjacent to all but $i+2$ of the colors of $X$,
so by  $Q(i)$ the at most $(i+2)$-coloring of $\partial (b) \cong S^i$ may be extended to cover $b$
 with $(i+2)$-colored  $(i+1)$-cells, each with diameter $< \epsilon$.  
Repeating for all $(i+1)$-cells of $X$ covers $X^{(i+1)}$.  
Lastly for each region $r \in X^{\langle n \rangle}$, 
$\partial (r) \cong S^{n-1}$ so by$Q(n-1)$ $r$ may be covered by small $(n+1)$-colored $n$-cells.  
Repeating for all the regions gives the $(n+2)$-colored covering of $S^n$.
\end{proof}
\begin{proof} (Theorem~\ref{small_cells}) (1) Let $\epsilon > 0$, let $r$ be one of the regions of $X$  having diameter $> \epsilon$, 
and assume there are $i$ of these larger regions. 
There is a tubular neighborhood $U=N_\delta(r) \cong B^n$ with $S^{n-1}=\partial B$ fully transverse to $X$.
$S^{n-1}$ inherits a local $(n+1)$-coloring from $X$.  Cutting along $S^{n-1}$, and using
 $Q(n-1)$ 
of Lemma~\ref{PQ} to fill in $\mbox{Int}(U)$,  gives $X'$ 
obtained by a surgery on $B^n$ that covers and replaces $r$ with
 only small regions.
Since $X'$ has one fewer regions with diameter $ < \epsilon$  than $X$ does,
and $X'\approx X$,
by induction this proves the claim.
For (2),
the tunneling step
 in Proposition~\ref{exists_cpx}
is a holonomy surgery, where $B^n$ is an $\epsilon$-neighborhood 
of the regions $r_1, r_2$ and their connecting edge in the proposition.
\end{proof}
\begin{thm}\label{make_to_ab}
Let $(M^n,N^{n-1})$ with possible boundary, $X$ ecpx on $M$, and 
non-separating simply connected fully transverse hypersurface $L^{n-1} \subset M$
be given.
Let generating set $\langle \gamma_1,\dots ,\gamma_g \rangle = \pi_1(M)$ be given with 
$\gamma_1 \cup L = \mbox{ pt.}$ and  $\gamma_i \cup L = \emptyset, i>1$. 
Let $\rho \in S_{n+1}$ be a permutation.
Then $\exists \bar{X}$ ecpx on $M$ with 
$\mbox{Hol}_{\bar{X}} (\gamma_1)  = \rho$
and 
$\mbox{Hol}_{\bar{X}}  (\gamma_i)  = \mbox{Hol}_X  (\gamma_i), \  \forall i \neq 1$.
\end{thm}

\begin{proof}
Color around the base point and 
let $a,b \leq n+1$ be distinct colors, and $(ab) \in S_{n+1}$ denote the transposition.  
We will first show
$\exists \bar{X}$ ecpx on $M$ with $\mbox{Hol}_{\bar{X}} (\gamma_1)  = (ab) \circ  \mbox{Hol}_X (\gamma_1) $.
Smooth $L$ to have tubular neighborhood $L_\delta$ and choose $X' \approx X$
with cell diameter $< \delta / 10$.
Since $\pi_1(L) = 1$ the regions  in $L_\delta$ may be   $(n+1)$-colored.
Choose this coloring $f$ by following along $\gamma_1$ to agree with the base point.
$X'$ may be isotopied keeping cell size $<\delta/5$
to have its $a$-colored regions missing $L$, 
for example by shrinking the $a$-colored regions to small neighborhoods of points.

Assume that $X'$ may be isotopied (keeping the $a$-colored regions away from $L$)
so that for each $b$-colored region $r$ of $X'$, $r \cap L$ 
is a collection of disks $D^{n-1}$ with $\partial D = S^{n-2} \subset \partial r$.
Since $L$ is simply connected it separates $L_\delta$ into distinct sides 
$L_1, L_2$ with $L_1 \cup L \cup L_2 = L_\delta$.

Let $\bar{X}$ be formed from $X'$ by inserting the disks $r \cap L$ for each $b$ colored region $r$ of $X'$.
Since these are disks and meet $X'$ transversely, $\bar{X}$ is a cpx on $M$.
We claim $\bar{X} \cap L_\delta$ may be $(n+1)$-colored as follows.

  Let $\bar{f}: \bar{X}^{\langle n \rangle} |_{L_\delta} \to \{1,\ldots,n+1\}$ be the 
local coloring of $L_\delta$ defined as follows:

\vspace{10pt}

\hspace{5em}
\parbox{3.5cm}{$\bar{f}(r) =f(r)\\
                                               \bar{f}(r) =f(r)\\   
                                               \bar{f}(r) =  \{a,b\} \setminus f(r)$} 
\hspace{20pt}
\parbox{5cm}{if $f(r) \not\in \{a, b \} $ \\
                         if $f(r) \in \{a, b \}$ and $r \subset L_1 $  \\
                         if $f(r) \in \{a, b \}$ and $r \subset L_2 $ }    
\hfill

\vspace{10pt}
\noindent   
Thus $\bar{f}$ is constructed from $f'$ by switching 
the $\{a,b\}$-colored regions on the $L_2$ side.  
One checks easily $\bar{f}$ is a coloring.                  
Thus   $\bar{X}$ is an ecpx on $M$.
Also 
$\mbox{Hol}_{\bar X}(\gamma_1) = (ab) \circ \mbox{Hol}_X(\gamma_1)$
and
$\mbox{Hol}_{\bar X}(\gamma_i) = \mbox{Hol}_X(\gamma_i), \ i>1$.

For $L$ transverse to $X'$ but meeting the $b$-colored regions in non-disk components,
the procedure may be followed with the $(n+1)$-coloring of 
the regions of $\bar{X}$ but with the exception that the subdivided regions
may not be homeomorphic to $B^n$.
Local surgeries along $S^{n-1} \subset L_\delta$ containing the deviant regions can 
replace $\bar{X}$ with a $(n+1)$-colored ecpx on $M$
 whose cells agree with $X'$ near $\partial L_\delta$.
The effect on $\gamma_1$ remains to interchange the $a$ and $b$ colors, and leave the other $\gamma_i$ unaffected.

By using multiple parallel copies of $L$ in Theorem~\ref{make_to_ab},
and repeating over $(ab) \circ (a'b') \circ \cdots \circ (a''b'') = \rho$,
 $\bar{X}$ can be constructed so Hol$_{\bar{X}}(\gamma) = \rho$
\end{proof}

Note that the proof of Theorem~\ref{make_to_ab} goes through for
$L$ not simply connected if
$\mbox{Hol}_X(\pi_1(L))=1$.

\begin{defi}
Let $X$ ecpx on $(M^n,N^{n-1})$.  $X$ is rigid if 
$\pi_1(X^{(1)}|_{ M \setminus N}) \to \pi_1(M)$ is surjective.
\end{defi}

\begin{thm}\label{make_hol_f}
Let $M^n$ be given with $f:\pi_1(M) \to S_{n+1}$.
Then $\exists X$ ecpx on $M$ with $Hol_X = f$.
\end{thm}

\begin{proof}
Let $f:\pi_1(M) \to S_{n+1}$ be given.  Choose $Z$ cpx on $M$ and $\delta >0$ s.t. 
$Z^{(1)}_\delta$ is a tubular neighborhood.
Choose $W$ ecpx on $M$ with cell diameter $<\delta/10$,
and put $W_1 = \{r \in W^{\langle n \rangle}| r \subset Z^{(1)}_\delta \}$.
Do this so $Z^{(1)}$ is a retract of $W_1$.  Note $W_1$ is rigid.

Choose base point $v \in Z^{(0)}$ and tree $T \subset Z^{(1)}$.
Let $\{e_1, \dots e_k\} = Z^{\langle 1 \rangle} \setminus T$.
  There is a 
$\gamma_1 \in \pi_1(T \cup e_1)$ that crosses $e_1$ once.  If 
$\mbox{Hol}_{Z^{(1)}}(\gamma_1) \neq f(\gamma_1)$, 
choose hyper-disk normal to $e_1$ that cuts 
$W_1 \cap N_\delta(e_1)$.  Apply Theorem~\ref{make_to_ab} 
to get $\bar{W_1}$,
with $\mbox{Hol}_{\bar{W_1}}(\gamma_1) = f(\gamma_1)$.
Repeat inductively for each of the remaining edges of $\{e_1, \dots e_k\}$,
choosing $\gamma_i$ to cross $e_i$ once 
to get $\hat{W}_1$, 
with 
$\mbox{Hol}_{\hat{W}_1}(\gamma_i)=f(\gamma_i)$.

Next let $b \in Z^{\langle 2 \rangle}$ be given.
Since $\pi_1(M)$ and $\pi_1( \hat{W}_1)$ are generated by $\{\gamma_1, \dots \gamma_k\}$,
and
$f$ and $\mbox{Hol}_{\hat{W}_1}$ agree on  $\{\gamma_1, \dots \gamma_k\}$, \
$\pi_1(\partial b) =1 \Rightarrow f(\partial b)=1  \Rightarrow \mbox{Hol}_{\hat{W}_1}(\partial b)=1$.
Thus the regions of  $\hat{W}_1$ meeting $\partial b$ can be 
 $(n+1)$-colored.  By Lemma~\ref{PQ}
the coloring may be extended to regions that cover $b$.  Repeating this for all $2$-cells yields $W_2$.
By induction on the cell dimension we get $W_n$ ecpx on $M$ with $\mbox{Hol}_{W_n}=f$.
\end{proof} 
The next two lemmas are needed to prove Theorem~\ref{hol_hol}
\begin{lem}\label{B^n_cut_fill}
Let $X, Y$ ecpx on $M^n$.  Let $B^n$ and $B_\epsilon = N_\epsilon(B)$ be given s.t. 
$\partial B$ and $\partial B_\epsilon$
are fully transverse to $X$ and $Y$.  Then there are $\bar{X}$ and $\bar{Y}$ ecpx on $M$ with
 $\bar{X} \approx X$, $\bar{Y} \approx Y$ and  $\bar{X}|_B = \bar{Y}|_B$.
\end{lem}
\begin{proof}
Choose ecpx $W$ on $B$.  $(n+1)$-color $W$.  Consider $(X \setminus B_\epsilon) \cup W$.
Since $\pi_1(\partial B_\epsilon)=1$ we may $(n+1)$-color the faces of the regions of $X \setminus B_\epsilon$
exposed to $B_\epsilon$.  By Proposition~\ref{fill_from_bdry}, the cavity $B_\epsilon \setminus B$ may be filled with 
$(n+1)$-colored regions extending $X \setminus B_\epsilon \cup W$ to form  $\bar{X}$.  
Then  $\bar{X} \approx X$.
Likewise choose $\bar{Y} \approx Y$ with $Y|_B = W$.
\end{proof}

\begin{lem}\label{disk_cut_fill}
Let $N^n \subset M^n$,
$N$ rigid and connected,
 and $X, Y$ ecpx on $M$ with $X|_N = Y|_N$.
Let $D^k \subset M \setminus \mbox{Int}(N)$ be a disk with 
$\partial D^k = S^{k-1} \subset \partial N$.
Then if $k\geq 2$ or $\mbox{Hol}_X =\mbox{Hol}_Y$,
there are $\bar{X} \approx X$, $\bar{Y} \approx Y$ s.t.
$\bar{X}$, $\bar{Y}$ agree on $N \cup N_\epsilon(D)$.
\end{lem}
\begin{proof}
Keeping $X|_N$ and $Y|_N$ fixed, we may find holonomic replacements $X'$ and $Y'$ 
with small cell diameter in $M \setminus N$.
Let $D_\epsilon = D \times [-\epsilon, \epsilon]$ 
and $D_{2\epsilon}=  D \times [-2\epsilon,2 \epsilon]$ 
denote the tubular neighborhoods of $D$ that only intersect $N$ on their ends.
Again by holonomy 
keeping 
$X'|_N$ and $Y'|_N$ fixed,
we may find $X''$ with $D^k \subset X''^{(k)}$
so that $X'' \setminus D_{2\epsilon}$ is a cpx, 
and find $Y''$ fully transverse to $D$ with $Y'' \cap  D_\epsilon$  a cpx.
Form $X'''$ by taking
$(X'' \setminus D_{2\epsilon}) \cup (Y'' \cap  D_{\epsilon})$.
Note $X'''$ has a cavity $ D_{2\epsilon} \setminus D_{\epsilon}$.
If $k\geq 2$ then $\pi_1(\partial X''') = 1$
and $\mbox{Hol}_{\partial X'''} = 1$ so 
$\partial X'''$ may be $(n+1)$-colored.
If $k = 1$ put $\partial D^1 = S^0 = \{a,b\}$.
Since $\mbox{Hol}_X =\mbox{Hol}_Y$,
 $X|_N = Y|_N$,
and $N$ is connected,
paths from $a$ to $b$ along $D^1$ or on $\partial D_{2\epsilon}$
will give the same coloring assignment.
So 
$\mbox{Hol}_{\partial X'''} = 1$ and
$\partial X'''$ may be $(n+1)$-colored.

Next, by Proposition~\ref{fill_from_bdry}, we may extend $X'''$ 
by filling the cavity with $(n+1)$-colored regions.  
This gives $\bar{X}$ ecpx on $M$,
with $\bar{X} \approx X'''$.
Also $\bar{X}|_{N \cup D_\epsilon} = Y''|_{N \cup D_\epsilon}$
and $Y'' \approx Y$.
\end{proof}

\begin{thm}\label{hol_hol}
Let $M^n$ be given with $X,Y$ ecpx on $M$.
Then  $Hol_X = Hol_Y \iff X \approx Y$
\end{thm}
\begin{proof}
In Proposition~\ref{approx_to_Hol} we showed $\Leftarrow$.
Let $X,Y$ ecpx on $M$ with $\mbox{Hol}_X = \mbox{Hol}_Y$.
Let $W$ ecpx on $M$.
Choose $X' \approx X$ and $Y' \approx Y$ so that $W$ has considerably larger cells 
than $X'$ and $Y'$
and $W$ is fully transverse to $X'$ and $Y'$.
Choose edge tree $T$ on $W^{(1)}$, and neighborhood $T_{2\epsilon}$ 
so that $T_{2\epsilon} \cong B^n$ and $\partial T_{2\epsilon}$ is fully transverse to $X'$ and $Y'$.
By Lemma~\ref{B^n_cut_fill} we may choose 
$X'' \approx X'$ and $Y'' \approx Y'$ so that 
$X''|_{T_\epsilon} = Y''|_{T_\epsilon}$.

Let $e \in W^{(1)}\setminus T$.  By Lemma~\ref{disk_cut_fill},
$\exists X''' \approx X''$ and $Y''' \approx Y''$ so that for some $\epsilon > 0$,
$X''' $ and $Y''' $ agree on $T_\epsilon \cup e_\epsilon$ .  By induction we may find
$\hat{X} \approx X''' \approx X$ and $\hat{Y} \approx Y''' \approx Y$ that agree on $N_\epsilon(W^{(1)})$
for some $\epsilon > 0$.

Let $i\geq 1$ and assume we have found 
$\hat{X} \approx X$ and $\hat{Y} \approx Y$ that agree on $N_\epsilon(W^{(i)})$. 
Let $b \in W^{(i+1)}$.  By Lemma~\ref{disk_cut_fill} we may find
$\hat{\hat{X}} \approx \hat{X}$ and $\hat{\hat{Y}} \approx \hat{Y}$ that agree on $N_\epsilon(W^{(i)}) \cup b_\epsilon$.
By induction on the number of $(i+1)$-cells, this can be extended to  $N_\epsilon(W^{(i+1)})$.
By induction on $i\leq n$ we have $\tilde{X}, \tilde{Y}$ ecpx on $W^{(n)} = M$
with $\tilde{X} \approx X$ and  $\tilde{Y} \approx Y$ s.t. $\tilde{X} = \tilde{Y} $.

So $X \approx Y$.
\end{proof}

\section{ Multi-layered cpx on $M^n$}       

We next investigate when several cpx are overlaid transversely on $M^n$.
\begin{defi}
Let ${\bf X}=(X_1, \dots ,X_j)$ be a family of setwise transverse cpx   on $M$.
Let ${\bf X}^{\langle n \rangle} = \bigcup_{i\leq j} X^{\langle n \rangle}_i$.
If $\forall Q \subset {\bf X}^{\langle n \rangle},  \ \bigcap Q \neq  \emptyset \Rightarrow \bigcap Q \cong B^i$, for some $i \leq n$, 
then say 
$(X_1, \dots , X_j)$ are fully transverse and
${\bf X}$ is a $j$-layer cpx on $M$.
\end{defi}

The definition says that any set of regions with non-empty intersection must intersect in an $i$-cell, for some $i$.
In particular, two regions on the same layer must be disjoint or have connected intersection.
A cpx $X$ is fully transverse as a $1$-layer cpx if there is a bijection between the cells of $X$ and the subsets of $X^{\langle n \rangle}$
with non-empty intersection.  Thus the definition is equivalent to requiring that for any subset of cells of arbitrary dimension in ${\bf X}$,
with at most one from each layer,
their intersection must be
in the form $B^i$ or $\emptyset$.

\begin{defi}
Let ${\bf X}=(X_1, \dots ,X_j)$ be a $j$-layer cpx on $M^n$.  
A $c$-coloring of ${\bf X}$ is a map
$f: {\bf X}^{\langle n \rangle} \to \{1,\ldots, c\} $  s.t.
$\forall \mbox{ distinct } r, r' \in {\bf X}^{\langle n \rangle}, \ r \cap  r' \neq \emptyset \Rightarrow   f(r) \neq  f(r')$.
Then say  ${\bf X}$ is $c$-colorable.  
If $\forall r \in {\bf X}^{\langle n \rangle}, \ {\bf X}|_{N_\epsilon(r)}$
 is $c$-colorable then say ${\bf X}$ is locally $c$-colorable.
\end{defi}
So adjacent regions in the same layer, and overlapping regions in different layers, must be given distinct colors.
It may be, if each region of $X_i$ intersects each region of $X_{\hat{i}}$ for example, that
the chromatic number of ${\bf X}$ is the sum of the chromatic numbers of the $X_i$.

\begin{prop}
Let $M^n$ and $j\geq 1$ be given.  Then $\exists$ a $j$-layer cpx ${\bf X}$ on $M^n$ with a 
$(n+j)$-coloring.  Moreover the $X_i$ may be chosen so each $X_i$ requires only $n+1$ colors in the coloring.
\end{prop}
 
\begin{proof}
Choose fully transverse ecpx $X_1$ on $M$ with a $(n+1)$-coloring $f:X_1^{\langle n \rangle} \to \{j, \dots, n+j \}$.
Choose $\epsilon >0$ so that the region faces may be expanded or contracted by $j \cdot \epsilon$
without changing the cell structure.
Construct $X_2$ offset from $X_1$, by expanding the $j$-color regions of $X_1$ by $\epsilon$,
then expanding the $(j+1)$-color regions of $X_1$ by $\epsilon$ into the remaining space, 
and continuing inductively until the $(n+j)$-colored regions of $X_2$ are slightly smaller
in every direction then the corresponding regions of $X_1$.
Now color each region of $X_2$ with $f:X_2^{\langle n \rangle} \to \{j-1, \dots, n+j-1 \}$, 
where the new regions are assigned $1$ less than the regions of $X_1$.
One checks this is an $(n+2)$-coloring of $(X_1, X_2)$ and the cpx are fully transverse.

Applying the procedure in the same manner to make $X_3$ from $X_2$ yields 
 coloring
$f:X_3^{\langle n \rangle} \to \{j-2, \dots, n+j-2 \}$, and an $(n+3)$-coloring 
on the fully transverse family  $(X_1, X_2, X_3)$.
By induction one gets a fully transverse family
 $(X_1, \ldots X_j)$ with $(n+j)$-coloring, where for each layer $X_i$
the regions use colors  $\{j-i+1, \dots, n+j-i+1 \}$.
\end{proof}

We will show each $j$-layer cpx on $M^n$ correspond to a cpx on $M^{n+j-1}$.
Define $K=K(B^{j-1})$ to be the dual to the standard $(j-1)$-simplex on $B^{j-1}$.
$K$ has one vertex $v_0 \in K^{(0)} \cap \mbox{Int}(B)$, and 
we may choose bijective coloring $f: K^{\langle j-1 \rangle} \to \{1, \ldots , j \}$.
Moreover there is a bijection $\gamma$ between the non-empty 
$J \subset \{1, \ldots , j\}$
and the various dimension cells $b$ in $K$ with $b \not\subset \partial B$, 
given by
$\gamma(J) = \bigcap f^{-1} (J) = b$.
$\gamma$ will be used in Theorem~\ref{main_layer}.

Let $M^n$ be a closed manifold, and consider the submanifold structure $M \times K(B^{j-1})$.
It may be further subdivided into a cpx $Y$ on $M \times B^{j-1}$.
Let $v_0$ be the interior vertex of $K$. Then we may choose $\epsilon > 0$ s.t.
$\forall \delta < \epsilon, \ Y|_{N_\epsilon (M \times v_0)} \cong  Y|_{N_\delta (M \times v_0)}$.
A $j$-layer cpx on $M^n$ corresponds to a cpx like $ Y|_{N_\epsilon (M \times v_0)}$
with the added condition that all cells lie either in the $M \times K$ submanifold structure
or are orthogonal to it.

\begin{thm}\label{main_layer} 
Let ${\bf X}=(X_1, \dots ,X_j)$ be a $j$-layer cpx on $M^n$. 
Identify $M=M\times \{v_0\} \subset M \times B^{j-1}$.
For $q \in {\bf X}^{\langle n \rangle}$ define
$\mbox{ind}(q) = i$ where $q \in X_i^{\langle n \rangle}$.
For each $Q \subset {\bf X}^{\langle n \rangle}$ with $\bigcap Q \neq \emptyset$,  
define $\Gamma_Q = \bigcap Q \times \gamma(\mbox{ind}(Q))$.
Then 
\[Y = \{\Gamma_Q\big| Q \subset {\bf X}^{\langle n \rangle}, \ \cap Q \neq \emptyset\} \
\cup  \ \{\Gamma_Q|_{M \times \partial(B^{j-1})}\big| Q \subset {\bf X}^{\langle n \rangle}, \cap Q \neq \emptyset\} \]
is a cpx on $M \times B^{j-1}$, with
$\Gamma_Q \in Y^{\langle n+j-|Q| \rangle}$.
\end{thm}

\begin{proof}
 Let $Q \subset {\bf X}^{\langle n \rangle}$ with $\cap Q \neq \emptyset$.  
Since the non-empty intersection of regions in layered cpx are $i$-cells, for some $i$,
and $\gamma$ produces  $i$-cells, for some $i$, 
the product  $\Gamma_Q = \cap Q \times \gamma(\mbox{ind}(Q))$
is homeomorphic to $B^i$, for some $i \geq 0$.
Write $Q=Q_1 \cup Q_2 \cup \cdots Q_j$ where $Q_i \subset X_i^{\langle n \rangle}$
and put $J = \{i|Q_i \neq \emptyset\}$.
By transversality in $M$, 
\[\mbox{dim}(\cap Q) = n - \sum_{i\in J} (|Q_i|-1) = n+|J| - |Q| \] 
Also
\[\mbox{dim}(\gamma(\mbox{ind}(Q))) = \mbox{dim}(\gamma(J))  = j - |J|\]
So  
\[\mbox{dim}(\Gamma_Q)=  \mbox{dim}(\cap Q)   +   \mbox{dim}(\gamma(\mbox{ind}(Q)))     = n+j-|Q| \]

\noindent
In particular $|Q| \leq n+j$.
For the vertices of $Y$ one may construct more, as  $\mbox{dim}(\Gamma_Q)=0 \Rightarrow 
\mbox{dim}(\gamma(\mbox{ind}(Q)))=0
\Rightarrow |J|=j$, and thus 
$\Gamma_Q = \cap Q \times \gamma(\mbox{ind}(1, \ldots j)) =
\{\mbox{pt}\} \times v_0 \in M \times K$.
This constructs $Y^{(0)} = \{ \Gamma_Q \big| \cap Q \neq \emptyset  \wedge |Q| = n+j\} \subset M \times v_0$.

Let $\Gamma_Q $ be a vertex and $q \in Q$.  
One checks $\Gamma_{Q \setminus \{q\}}$ is an edge with endpoint $\Gamma_Q$,
where if $\mbox{ind}(Q) = \mbox{ind}(Q\setminus \{q\})$ the edge is in the $M$ direction,
while if  
 $\mbox{ind}(q) \not\in \mbox{ind}(Q  \setminus \{q\} )$, then
$\gamma(\mbox{ind}(Q  \setminus \{q\}  ))$ is a $1$-cell of $K$ and 
$\Gamma_{Q  \setminus \{q\} } = \{pt\} \times \gamma(\mbox{ind}(Q  \setminus \{q\} ))$. 
Since there are $(n+j)$ choices for $q \in Q$, this shows the degree of the vertices is $(n+j)$.

For induction let $Q \subset {\bf X}^{\langle n \rangle}$ with $\cap Q \neq \emptyset$ and $|Q| = m$.
One may show $\Gamma_Q$ is a $(n+j-m)$-cell, meeting $(m-1)$ of the $(n+j-m+1)$-cells,
and $Y$ has regular degree.  In particular the regions of $Y$ are generated by
singletons $Q \subset  {\bf X}^{\langle n \rangle}$,
where $\Gamma_Q = Q \times \gamma(\mbox{ind}(Q))$ and $\gamma(\mbox{ind}(Q))$ fills
$1/j$ of the volume cone angle in $B^{j-1}$.
\end{proof}

\begin{defi}
Define $\Gamma: \{ Q \subset {\bf X}^{\langle n \rangle} | \cap Q \neq \emptyset\} \to Y^{\langle n+j \rangle}$
as in Theorem~\ref{main_layer}, and let $\Gamma( {\bf X})$ denote $Y$.
\end{defi}

\begin{prop}\label{col_col}
Let ${\bf X}=(X_1, \dots ,X_j)$ be a $j$-layer cpx on $M^n$ and $Y = \Gamma({\bf X})$ the corresponding cpx on $M^n \times B^{j-1}$.
Then  ${\bf X}$ is $c$-colorable $\iff$ $Y$ is $c$-colorable.
\end{prop}
\begin{proof}
Define $F: {\bf X}^{\langle n \rangle}  \to Y^{ \langle n+j \rangle}$ by
$F(r)= \Gamma_{\{r\}} = r\times \gamma(\mbox{ind}(r))$.  Note $F$ is $\Gamma$ restricted to singletons and is bijective.
Moreover $\forall r, r' \in {\bf X}^{\langle n \rangle}, \ r \cap r' \neq \emptyset \iff F(r) \cap F(r') \neq \emptyset$.  
Thus the claim follows.  Moreover, $F$ gives a bijection between regions of the same color.
\end{proof}

\section{$S_{n+j}$ Holonomy on $j$-layered ecpx on $M^n$}   

The holonomy results for single layer ecpx on $M^n$ follow through to similar results 
for $j$-layer ecpx on $M^n$, either by applying a similar proof 
or by applying Theorem~\ref{main_layer} and Proposition~\ref{col_col}.  
The results are restated below and proofs sketched.

\begin{defi}  Let ${\bf X}$ be a $j$-layer cpx on $M^n$.
Define $\|{\bf X}\|$ to be the CW-complex on $M$ whose cells are of the form
$\cap R$, where $R \subset {\bf X}^{\langle n \rangle}$,
$\cap R \neq \emptyset$,
and $R$ contains at least one region from each layer of ${\bf X}$.
\end{defi}

Each intersection above is of the form $B^i$, for some $ i\geq 0$, 
and one may check this is the coarsest CW-complex that refines each of the $X_i$ in ${\bf X}$.
A cell $b= \cap R$ of $\|X\|$ satisfies the regular degree condition iff at most one layer of $R$ has multiple regions.
\begin{defi}
For ${\bf X}$ a $j$-layer cpx on $M^n$, say  ${\bf X}$ is a $j$-layer ecpx on $M$
 if all $2$-cells of $\|{\bf X}\|$ are even-sided.
\end{defi}
\begin{prop}
 Let ${\bf X}$ be a $j$-layer cpx on $M^n$, and $Y= \Gamma({\bf X})$.
Then ${\bf X}$ is a  $j$-layer ecpx on $M$ $\iff$ $Y$ is an ecpx on $M \times B^{j-1}$.
\end{prop}
\begin{proof} Let $R\subset {\bf X}^{\langle n \rangle}$, with $\cap R \neq \emptyset$.
There is only one cell of $K(B^{j-1})$ that does not meet $\partial B^{j-1}$,
namely $v_0 \in K(B^{j-1})$
and moreover $\gamma(ind(R))= v_0 \iff R$ contains a cell from each layer.
Thus 
$\cap R$ is a cell of $\|{\bf X}\|$ 
$\iff$
$\mbox{ind}(R) = \{1,\dots , j\}$
$\iff$
$\Gamma(R) = \cap R \times \{v_0\}$
$\iff$
$\Gamma(R) \cap (M \times \partial B) = \emptyset$.
Thus all $2$-cells of $\|{\bf X}\|$ are even
$\iff$ all $2$-cells of $Y$ not meeting $M \times \partial B$ are even
and the claim follows.
\end{proof}

As in Definition~\ref{defi_HolW}, we may define for  ${\bf X}$ a $j$-layer cpx on $M^n$,
and  $ W$ a tubular neighborhood of $\|{\bf X}\|^{(1)}$, 
a holonomy map 
\[ \mbox{Hol}_W: \pi_1(W) \to S_{n+j}  \]
At each point on $M$, ${\bf X}$ has $j$-layers. 
Paths in $W$ carry  the $n$-colors surrounding an edge of a given layer
along with the $j-1$-colors of the other layers.  More generally, since each edge of $\|{\bf X}\|^{(1)}$
is the intersection of $n+j-1$ regions of ${\bf X}$,
there is one remaining color to apply at the ends of an edge and the holonomy is forced as  in the one-layer case.
The proof can be shown by considering $\Gamma({\bf X})$.

\begin{prop}\label{j_HolX}
For ${\bf X}$ a $j$-layer ecpx on $M^n$ 
the holonomy extends uniquely from
 $\mbox{Hol}_{{\bf X}^{(1)}}: \pi_1({\bf X}^{(1)}) \to S_{n+j}$ to
\[ \mbox{Hol}_{\bf X}: \pi_1(M) \to S_{n+j}\]
\end{prop} 
\begin{proof}
As the correspondence ${\bf X} \leftrightarrow \Gamma({\bf X})$ preserves local $(n+j)$-coloring the claim follows.
\end{proof}

\begin{thm}\label{j_loc123}
Let ${\bf X}$ be a $j$-layer cpx on $M^n$.  Then
\begin{enumerate}
\item        ${\bf X}$ is locally $(n+j)$-colorable $\iff$ ${\bf X}$ is an ecpx.
\item If $M$ is orientable then:
               ${\bf X}$ is locally $(n+j)$-colorable $\iff$  $\|{\bf X}\|^{(1)}$ is even cyclic.
\item If $M$ is simply connected then:
               ${\bf X}$ is locally $(n+j)$-colorable $\iff$  ${\bf X}$ is $(n+j)$-colorable.
\end{enumerate}
\end{thm}
\begin{proof}
The claims are similar to Theorem~\ref{loc123}.
(1) and (3) follow from Proposition~\ref{j_HolX}.  
For (2), first
${\bf X}$ is orientable $\iff$ $\Gamma({\bf X})$ is orientable.
Then
${\bf X}$ is locally $(n+j)$-colorable $\iff$ 
$\Gamma({\bf X})$ is locally $(n+j)$-colorable $\iff$
$\Gamma({\bf X})^{(1)}$ is even cyclic $\iff$
 $\|{\bf X}\|^{(1)}$ is even cyclic and the claim follows.
\end{proof}

\begin{prop}\label{j_fill_from_bdry} 
Let $(M^{n},N^{n-1})$ be given with  ${\bf Y}$ an $(n+j)$-colorable $j$-layer cpx on $N$.
Then ${\bf Y}$ has an extension to an $(n+j)$-colorable ecpx ${\bf X}$ on $M$.
\end{prop}
The general argument follows from Proposition~\ref{fill_from_bdry}
with exception that Lemma~\ref{Sn_fill} is reworked below.

\begin{lem}\label{j_Sn_fill}
Let ${\bf X}$ be a $j$-layer cpx on $S^n$ with a $(n+j+1)$-coloring.  
Then ${\bf X}$ extends to a $j$-layer cpx on $B^n$ with $(n+j+1)$-coloring.
\end{lem}
\begin{proof}
Let ${\bf X}$ be a $j$-layer cpx on $S^n$ with $(n+j+1)$-coloring.  
For initial step $n=0$, let $S^0 =\{ u, v\}$ and ${\bf X}= \{u_1, \dots , u_j, v_1, \dots v_j\}$
have corloring $f:{\bf X} \to (1,2, \ldots,j+1)$.
As seen directly or by Lemma~\ref{Sn_fill}, the first layer $ X_1 = \{u_1,v_1\}$
may be extended to $Y_1$ cpx on $B^1$.
Assume for some $i$, $0\leq i <  j$, that  $(X_1,\ldots, X_i)$ has been extended to an $i$-layer
$(n+i+1)$-colored
cpx $(Y_1, \ldots, Y_i)$ on $B^1$.
For each $w \in \|(Y_1, \dots Y_i)\|^{\langle 0 \rangle}$ in Int$(B^1)$,
for some $\epsilon >0$, $N_\epsilon(w)$ meets $(i+1)$ regions in $(Y_1, \dots Y_i)$.
Since there are $j+1 > i+1$ colors, we may begin to construct $Y_{i+1}$
by choosing a color for $N_\epsilon(w)$.  
Do this for each $w \in  \|(Y_1, \dots Y_i)\|^{\langle 0 \rangle}$ in Int$(B^1)$.
$Y_{i+1}$ is a series of colored dots, with connecting segments adjacent to $i$ regions in 
 $(Y_1, \dots Y_i)$.
This leaves $j+1 - i \geq 2$ colors for covering each connecting segment of $Y_{i+1}$.
By Lemma~\ref{Sn_fill}, $X_{i+1}$ may be extended to $Y_{i+1}$ on $B^1$
here with  $(Y_1, \dots Y_{i+1})$ an $(i+2)$-coloring cpx on $B^1$.

For induction on $n$, assume the condition holds for $(n+j+1)$-colored cpx on $S^n$ and let 
${\bf X}= (X_1, \ldots X_j)$ be an $(n+j+2)$-colored $j$-layer cpx on $S^{n+1}$.
By Lemma~\ref{Sn_fill}, $ X_1$ may be extended to an
 $(n+j+2)$-colored  cpx $ Y_1$ on $B^{n+1}$.
For induction on the $i$ layers,
Suppose $(Y_1, \dots, Y_i)$ is an  $(n+j+2)$-colored $i$-layer cpx  on $B^{n+1}$
that extends $(X_1, \ldots, X_i)$.
$\forall w \in \|(Y_1, \ldots Y_i)\|^{\langle 0 \rangle}$,
$N_\epsilon(w)$ meets $(n+i+1)$ regions of $(Y_1, \ldots, Y_i)$.
Since $(n+j+1)-(n+i+1) \geq 1$, each $N_\epsilon (w)$
may be colored forming $Y_{i+1}^{\langle 0 \rangle}$.
By Lemma~\ref{Sn_fill}, we may inductively in $k$  cover 
 $Y_{i+1}^{\langle k \rangle}$, to construct  $Y_{i+1}$.
Here $(Y_1, \ldots, Y_{i+1})$ is an $(n+j+2)$-colored $(i+1)$-layered cpx on $B^{n+1}$ that extends $(X_1, \ldots, X_{i+1})$.
Induction on $i$ until $i=j$  followed by induction on $n$ shows the claim.
\end{proof}

For ${\bf X}$ a $j$-layer ecpx  on $M^n$ we similarly define a local surgery 
that does not impact the holonomy of ${\bf X}$. 
Let $B^n \subset M^n$ be given with boundary $S^{n-1} $ fully transverse to ${\bf X}$.
Replace $j$-layer ecpx ${\bf Y}=B^n \cap {\bf X}$ with a choice of $j$-layer ecpx ${\bf Y'}$ on $B^n$,
where it is required
$S^{n-1}$ is fully transverse to ${\bf Y'}$,  $S^{n-1}  \cap {\bf Y'} = S^{n-1} \cap {\bf Y}$, and
${\bf Y'} \cup_{S^{n-1}} {\bf Y}$ is a $j$-layer ecpx on $S^n$.
Then respecting layers and
keeping $S^{n-1}$ fixed, we may cut out ${\bf Y}$ and sew in ${\bf Y'}$.
Note each cell of ${\bf Y}$ cut by $S^{n-1}$ is re-spliced with a partial cell of ${\bf Y'}$.
Call this a holonomy surgery on $B^n$.

\begin{defi}\label{j_Hol_Xi}
Let $(M^n, N^{n-1})$ be a manifold with possible boundary  $N$.  Let ${\bf X}$ be a $j$-layer ecpx and ${\bf X'}$ a $j$-layer cpx on $(M,N)$
that agree on an $\epsilon$ neighborhood of $N$.
If
there exists a finite number $k$ of $j$-layer cpx's ${\bf X_i}$ on $M$ with ${\bf X_0} = {\bf X}$ and ${\bf X_k} ={\bf X'}$,
and
with ${\bf X_i} \cap N_\epsilon = {\bf X} \cap N_\epsilon,  \forall i \leq k$
s.t.  $\forall i$ either
\begin{enumerate}
\item ${\bf X_i}$ and ${\bf X_{i+1}}$ are ambient isotopic keeping $N_\epsilon$ fixed.
\item ${\bf X_{i+1}}$ is obtained from a holonomy surgery  on ${\bf X_i}$  keeping $N_\epsilon$ fixed.
\end{enumerate}
Then write 
${\bf X} \approx {\bf X'}$, and define $[{\bf X}]$ to be the collection of all   $j$-layer cpx ${\bf X'}$ s.t. ${\bf X} \approx {\bf X'}$.
\end{defi}

\begin{prop}\label{j_approx_to_Hol} Let ${\bf X}$ be a $j$-layer ecpx on $M^n$ and ${\bf Y}$ a $j$-layer cpx on $M$.
 Then ${\bf X} \approx {\bf Y}$ implies  ${\bf Y}$ is an ecpx on $M$ and moreover
$\mbox{Hol}_{\bf X} =  \mbox{Hol}_{\bf Y}$.
\end{prop}
\begin{proof} The proof is similar to Proposition~\ref{approx_to_Hol}
\end{proof}

\begin{thm}\label{j_make_to_ab}
Let $(M^n,N^{n-1})$ with possible boundary, ${\bf X}$ $j$-layer ecpx on $M$, and 
non-separating simply connected fully transverse hypersurface $L^{n-1} \subset M$
be given.
Let generating set $\langle \gamma_1,\dots ,\gamma_g \rangle = \pi_1(M)$ be given with 
$\gamma_1 \cup L = \mbox{ pt.}$ and  $\gamma_i \cup L = \emptyset, i>1$. 
Let $\rho \in S_{n+1}$ be a permutation.
Then $\exists {\bf \bar{X}}$ $j$-layer ecpx on $M$ with 
$\mbox{Hol}_{\bf{\bar X}} (\gamma_1)  = \rho$
and 
$\mbox{Hol}_{\bf{\bar X}}  (\gamma_i)  = \mbox{Hol}_{\bf X}  (\gamma_i), \  \forall i \neq 1$.
\end{thm}
\begin{proof} The proof is similar to Theorem~\ref{make_to_ab}
\end{proof}

\begin{thm}\label{j_make_hol_f}
Let $M^n$ be given with $f:\pi_1(M) \to S_{n+1}$.
Then $\exists {\bf X}$ $j$-layer ecpx on $M$ with $Hol_{\bf X} = f$.
\end{thm}
\begin{proof} The proof is similar to Theorem~\ref{make_hol_f}
\end{proof}

\begin{thm}\label{j_hol_hol}
Let $M^n$ be given with ${\bf X},{\bf Y}$ $j$-layer ecpx on $M$.
Then  $Hol_{\bf X} = Hol_{\bf Y} \iff {\bf X} \approx {\bf Y}$
\end{thm}
\begin{proof} The proof follows Theorem~\ref{hol_hol}.
\end{proof}

\section{The 4-Color Theorem vs. the Poincar\'{e} Conjecture}\label{vs}    
This paper developed from interest in determining if there were similarities between
the  4-Color Theorem and the Poincar\'{e} Conjecture.  
One re-occurring goal was to see if the two famous theorems could be reformulated so that they were special cases of a general statement.
The propositions in this section show three attempts at comparing the two theorems in this way.

\begin{defi}\label{simp_ext}
Let $N^{n-1}$ be oriented, and $Y$ a cpx on $N$.  If $M^n$ is simply connected with boundary $N$,
and $X$ ecpx on $M$, with
 $X|_N = Y$, then say
$X$ is a simple extension of $Y$.
\end{defi}

\begin{prop}\label{infinity}
Consider the statements:
\begin{enumerate}
\item $\forall Y \mbox{ cpx on } S^2, \ \big| \{ [X] | X \mbox{ is a simple extension of } Y\} \big| > 0$
\item $\forall Y \mbox{ cpx on } S^2, \ \big| \{ [X] | X \mbox{ is a simple extension of } Y\} \big|  < \infty$
\end{enumerate}
Then it may be shown by short proofs that
\[ \mbox{ (1) } \iff \mbox{ the Four Color Theorem }\] 
\[ \mbox{ (2) } \iff \mbox{ the Poincar\'{e} Conjecture }\]
\end{prop}
\begin{proof}
Corollary~\ref{Sn_even} shows the $4$-Color Theorem  $\Rightarrow$  (1).
Let $Y$ be a cpx on $S^2$.
(1) $\Rightarrow \exists$ simple extension $X$ of $Y$.
By Theorem~\ref{loc123} $X$ is $4$-colorable, showing $Y$ is $4$-colorable.

For the second claim, starting with the Poincar\'{e} Conjecture\cite{PC},
we know $B^3$ is the only simply connected manifold with boundary $S^2$.
So (2) becomes the claim that each $Y$ cpx on $S^2$ has finitely many extension classes $[X]$ on $B^3$.
Each extension class corresponds to a coloring of $Y$, of which there are finitely many.  If however the Poincar\'{e} Conjecture were false,
and there 
was a simply connected $\tilde{S}^3 \not\cong S^3$, 
then choose $W$ ecpx on $\tilde{S}^3$.  
Removing a neighborhood of a vertex gives $W'$ ecpx on $\tilde{B}^3 \not\cong B^3$.  
$W'$ has a tetrahedral coloring on $\partial W' = S^2$.  Then choosing  a 4-colorable $Y$ cpx on $S^2$  
and simple extension $X$ on $B^3$,
$W'$ could be spliced in at any vertex inside $X$.  Doing this $i$ times produces the connected sum of $i$ copies of $\tilde{S}^3$.  
For distinct $i$ these are not homeomorphic and so not holonomic, showing infinitely many classes of simple extensions of $Y$.
\end{proof}

Note the second claim only needs the existence of a colorable cpx on $S^2$ with finitely many simple extensions.
The result can also be expresses using transverse extensions.
\begin{defi}\label{transv}
Let $N^{n-1}$ be oriented, and $Y$ a cpx on $N$. 
If $M^n$ is simply connected with boundary  $N$,
and $X$ ecpx on $M$, with
$X|_{N}$ transverse to  $Y$, 
and $X^{(1)} \cup Y^{(1)}$ even cyclic,
then say $X$ is a simple transverse extension of $Y$.
\end{defi}

In the definition the $1$-skeleton from all edges including those on the boundary of $M$ from $X$ and $Y$ are joined and then tested for evenness.
If
$M_\epsilon = M\cup_{ N\times \{0\}} (N\times [0,\epsilon])$
denotes the tubular neighborhood off the boundary by $\epsilon$,
then the definition is equivalent to saying
 $X \cup_{ N\times \{0\}} (Y\times [0,\epsilon])$ is an ecpx on $M_\epsilon$
or equivalently that it is a simple extension of $Y \times \{\epsilon\}$ from Definition~\ref{simp_ext}.
Proposition~\ref{infinity} can be rephrased as
\begin{cor}
The statements
\begin{enumerate}
\item $\forall W \mbox{ cpx on }S^2, \exists W' \mbox{ cpx on }S^2 \mbox{ s.t. }\\
\big|\{[X]|X \mbox{ is a simple transverse extension of } W \mbox{ and } X|_{S^2} = W'\} \big| \geq 1$
\item $\forall W \mbox{ cpx on }S^2, \forall W' \mbox{ cpx on }S^2 ,\\
\big|\{[X]|X \mbox{ is a simple transverse extension of } W \mbox{ and } X|_{S^2} = W'\} \big| \leq 1$

\end{enumerate}
are equivalent to the $4$-Color Theorem and Poincar\'{e} Conjecture, respectively.
\end{cor}

\begin{proof}
(1) Any $X$ meeting the set criterion gives a coloring.  Conversely any coloring can
be transversely extended to fill $B^3$.
(2) Suppose $W'$ and $X$ ecpx on $B^3$ satisfy the set criterion.
Then $(W,W')$ is a $2$-layer ecpx with forced $4$-coloring on $S^2$.
Suppose $\tilde{X}$ ecpx on $B^3$  with  $\tilde{X}$ and $W'$ also satisfying the set criterion.
Then $X$ and $\tilde{X}$ agree in the coloring on $S^2$ so $X \approx \tilde{X}$.
Thus Poincar\'{e} Conjecture $\Rightarrow$  $ \big| \{\ldots \} \big| \leq 1$.
If $ \exists$ simply connected $\tilde{S}^3 \not\cong S^3$,
then forming connected sums as in Proposition~\ref{infinity} would yield $ \big| \{\ldots \} \big| > 1$.
\end{proof}

The next comparison between the two famous theorems is helped by slightly enlarging the class of cpx on $M^n$ 
to include a trivial case.  Consider the $2$-coloring of $S^2$ into
southern and northern hemispheres.  This is a submanifold structure
consisting of the two hemispheres and the equator and satisfies the regular degree condition;
However it has no vertices.  For higher dimensional analogs
we will start constructing $X$ on $M^n$ similarly with $S^k$ as the $k$-skeleton,
no $j$-skeleton for $j<k$, 
and for $i\geq 1$ inductively form the $(k+i)$-skeleton by sewing on a set of 
$B^{k+i}$ to the $(k+i-1)$-skeleton.
If each cell is embedded in $M^n$ and satisfies the regular degree condition, then call 
$X$ a generalized coloring complex on $M^n$, and write $X$ gcpx on $M^n$.
While these are natural coloring objects it only serves to include structures on $S^n$:

\begin{lem}\label{No_vert}
Let $n\geq 1$ and $X$ be a gcpx on $M^n$,
\begin{enumerate}
\item If $X^{(0)}$ has two vertices $M^n \cong S^n$.
\item $X^{(0)}$ cannot have one vertex.
\item If $X^{(0)}$ has no vertices $M^n \cong S^n$.
\end{enumerate}
\end{lem}
\begin{proof}
(1) Let $\{v_0, v_1\} =X^{(0)}$ .
The regular degree requires $X$ to have $n+1$ edges, each with one endpoint at $v_0$ and one endpoint at $v_1$.
For $X$ smooth positioning points on these edges distance $\epsilon$ from $v_0$ and taking the convex hull 
gives a ball around $v_0$.  One can show there is an isotopy of the $S^{n-1}$ boundary of this ball
to an $\epsilon$ neighborhood of $v_1$ by isotopying the cells of the CW-complex $Y$  defining the sphere.
The points $Y^{(0)}$ can be isotopied along $X^{(1)}$  to $v_1$.
For $i \geq 1$, the isotopy of $Y^{(i-1)}$  can be extended to an isotopy of 
$i$-cells of $Y$ along $X^{(i+1)}$.
By induction this yields an isotopy of the sphere from a neighborhood of $v_0$ to a neighborhood of $v_1$
which demonstrates $M \cong S^n$.

For (2), if $X$ had only one vertex, no $1$-cells could be attached injectively, so $M=$ point.
For (3). In (1) we showed if  $X^{(0)} = S^0$ then $M^n \cong S^n$.  Let $0 \leq j \leq n$
 and suppose the claim has been shown for $X$ with minimal rank skeleton $X^{(j)} \cong S^j$
that $M^n \cong S^n$.
Let $X$ have minimal rank skeleton $X^{(j+1)} \cong S^{j+1}$.
As $X$ has regular degree there are $n-j$ regions adjacent to $ S^{j+1}$.
But every cell of rank $ > j+1$ is sewn onto  $ S^{j+1}$.
So $X$ has $n-j$ regions.
Choose one region $r$.  
We may choose a disk $D^{n-1} \subset r$ with $\partial D \subset \partial r$ and so that 
$D \cap  S^{j+1} = S^j$ and the intersection is  transverse.
$D$ bounds  two $n$-cells $r_1, \ r_2$ in $r$.
Taking $X' $
to be the gcpx on $M$ with cells generated by intersections of 
$X \setminus \{r\} \cup \{D, r_1, r_2\} $
and applying induction shows $M\cong S^{n}$.
\end{proof}

Note that while a CW-complex is usually constructed inductively from several $0$-cells,
there is no need to assume the minimal
$k$-skeleton consists of multiple $S^k$.  This follows since except for $B^1$, 
the boundary of $B^j$ is connected, and the sewing process cannot join multiple $S^k$
into the one connected component of $M^n$.

\begin{defi}If $M^n$ satisfies the property that
 for \underline{some} and \underline{all} locally $4$-colorable gcpx $X$ on $M$, $X$ is $4$-colorable,
then say $M$ is $4$-color extendable.
\end{defi}
\begin{prop}\label{Main_Sn}
For $n\geq 1$, let $M^n$ be closed connected with no boundary.  Then \\
For $n=1, \ 2, \mbox{ or } n\geq 4$,
\[M^n\cong S^n \iff M  \mbox{ is } 4\mbox{-color extendable} \]
For $n=3$,
\[M^n\cong S^n  \Rightarrow M  \mbox{ is } 4\mbox{-color extendable}  \Rightarrow H(M^n) = H(S^n) \]
\end{prop}
\begin{proof}
For $n=1$, $M^n = S^1$ and all $X$ are $4$-colorable.  For $n=2$, for every $M^2$
every cpx is locally $4$-colorable.  So the condition is simply that all $X$ on $M$ are $4$-colorable.
When $M=S^2$ this is the 4-Color Theorem.  When $M\neq S^2$,
it is well known that the complete graph $K_5$ embeds on $M$, 
and dual constructions yield cpx $X$ with chromatic number $\geq 5$.

For $n=3$, 
assume $M$ is simply connected.
Proposition~\ref{exists_cpx} states that there always exists a $4$-colorable cpx. 
Theorem~\ref{loc123} shows that when $M$ is simply connected, every locally $4$-colorable cpx
is in fact $4$-colorable.
This shows $M  \mbox{ is } 4\mbox{-color extendable}$.

For the second implication, assume $H(M) \neq H(S^3)$.
Then $H_1(M) \neq 0$
and we may choose a non-trivial map $H_1(M)\to S_2$.
Define $f:\pi_1(M)\to H_1(M) \to S_2 \hookrightarrow S_4$.
By Theorem~\ref{make_hol_f} there is an $X$ ecpx on $M$ with 
$\mbox{Hol}_X = f$.  Thus $X$ is not $4\mbox{-color extendable} $.

For $n \geq 4$, let cpx $X$ on $M^n$ be locally $4$-colorable.  
This implies $X$ can have no vertices
or other components of degree $ > 4$.  In fact, it has at most $4$ regions and so is $4$-colorable.
By Lemma~\ref{No_vert}, $M$ having no vertices implies $M\cong S^n$.
Conversely, if $M\cong S^n$, splitting $M$ into two hemispheres demonstrates the existence of a $4$-coloring.
\end{proof}
So for $n=2$ the proposition is equivalent to the $4$-Color Theorem, and for $n=3$ it shows $M$ is a homology sphere,
where homotopy sphere would be needed for $M^3 \cong S^3$
to be equivalent to the Poincar\'{e} Conjecture.

Instead of considering when local $4$-colorability implies $4$-colorability, the next definition
considers when forced local colorings extend on $j$-layer ecpx.

\begin{defi}
Let $M^n$ be given.  If every $j$-layer ecpx on $M$ is $(n+j)$-colorable, say 
$M$ is $(n,j)$-color extendable.
\end{defi}
Claim (3) below uses and is equivalent to the Poincar\'{e} Conjecture.
\begin{prop}\label{nj_extend}
Let $M^n$ be given.
\begin{enumerate}
\item If $ n \geq 2$.  Then  $ \pi_1(M^n)=1 \Rightarrow (\forall j \geq 1, M$
is $(n,j)$-color extendable) $\Rightarrow \pi_1(M)$ 
has no normal proper subgroup of finite index. 
\item If $n=2, j\geq 1$.  Then  
$ M^n \mbox{ is } (n,j)$-color extendable 
$\iff M=S^2  $ .
\item If $n=3$.  Then
$(\forall j \geq 1,  M^n \mbox{ is } (n,j)$-color extendable )
$\iff M=S^3 $     
\item If $n \geq 4$, \
$\exists M^n$ s.t. $\pi_1(M) \neq 1 \wedge \forall j>0, \ M$ is  $(n,j)$-color extendable.
\end{enumerate}
\end{prop} 
\begin{proof} (1) If $M$ is simply connected, then $\mbox{Hol}_{\bf X}$ is trivial and 
$M$ is  $(n,j)$-color extendable.  To show the second implication,
suppose $G \triangleleft  \pi_1(M)$ and $G$ is a proper subgroup of finite index $k$.  
Put $k'= \max(k,n+1)$.
Define
\[f: \pi_1(M) \to  \frac{\pi_1(M)}{ G} \to S_k \hookrightarrow S_{k'} \]
where the maps are the quotient, the group action, and the inclusion, respectively.
$f$ has non-trivial image.  Put $j=k' - n$.  
Thus by Theorem~\ref{j_make_hol_f} we may choose a $j$-layer ecpx ${\bf X}$ on $M$
with $\mbox{Hol}_{\bf X} = f$, showing $M$ is not  $(n,j)$-color extendable.

(2) Let $j \geq 1$  Except for $S^2$, all closed $2$-manifolds have non-trivial map
$f: \pi(M) \to H_1(M) \to \mathbb{Z}_2 \to S_{2+j}$.  
Thus $M\not \cong S^2 \Rightarrow M$ not  $(2,j)$-color extendable.

(3) For $M^3$ closed manifold, $\pi_1(M)$ is residually finite\cite{RF}.  A group $G$ is
 residually finite  if for every nontrivial element $g$ in $G$ there is a homomorphism $h$ from $G$ to a finite group, 
such that  $h(g) \neq 1$.   This means that either $\pi_1(M)=1$ or 
for some $ k\geq 2$,  $\exists \mbox{ non-trivial } f:\pi_1(M) \to S_k$.
But $M$  $(n,j)$-color extendable implies no such $f$ can exist.
Thus $\pi_1(M) = 1$.  The assertion that $M= S^3$ is thus true and equivalent to the Poincar\'{e} Conjecture.

(4) There are finitely presented simple groups of infinite order\cite{HG}, and choose such a group $G$.
For $4$-manifolds, $M^4$ may be chosen so $\pi_1(M)$ is any finitely presented group\cite{4M}.
Choose $M$ so that $\pi_1(M) = G$.
Then $\forall j \geq 1, \forall f: \pi_1(M) \to S_j$, $f$ must have a trivial image,
so $M$ is $(n,j)$-color extendable and not simply connected.
\end{proof}

The next attempted comparison while dependent upon coincidences is easily stated.

\begin{cor}
For $n\geq 2$, let $M^n$ be given.  
The following statement is true and equivalent to the 
$4$-Color Theorem and the Poincar\'{e} Conjecture
 in dimensions $2$ and $3$, respectively.
\[M\cong S^n \iff 
 M \mbox{ is } 4 \mbox{-color extendable and } M \mbox{ is } (n,j) \mbox{-color extendable}\]

\end{cor}
\begin{proof}
The hypothosis are combined from Propositions~\ref{Main_Sn} and \ref{nj_extend},
which are equivalent to the $4$-Color Theorem and  Poincar\'{e} Conjecture in dimension $2$ and $3$, respectively.
\end{proof}

\section{Applications}   

Two topological characterizations of the $4$-Color Theorem are repeated.

\begin{prop}Let $X$ be a cpx on $S^2$.  Then $X$ is $4$-colorable $\iff \exists X'$ cpx on $S^2$
s.t. $X$ and $X'$ are transverse and 
$\| (X , X') \|$ is even-cyclic 
$\iff \exists Y$ ecpx on $B^3$ that extends $X$.
\end{prop}
\begin{proof} This is a special case of Corollary~\ref{Sn_even} and Theorem~\ref{j_loc123}.
\end{proof}

For $X$ cpx on $S^3$ the chromatic number may be arbitrarily high.  
For example, $X$ may contain $k$-regions that are all adjacent to each other.
However when the problems with local colorability are dispersed, some bounds can be given.
 
Let $X$ be a cpx on $M^3$. Let $Q \subset X^{\langle 3 \rangle}$ be the set of regions $r$
with $N_\epsilon(r)$ not $4$-colorable.
Let $g$ be the graph that is subset of the dual of $X$ that has
a vertex centered in each region of $Q$
 and edges through the odd $2$-cells of $X$.
It is easily shown that $g$ has even degree vertices, since any $3$-cell must have  
an even number of odd sided $2$-cells.
Define $G\supset g$ to be the adjacency graph of $Q$.
$G$ has edges in additional to those of $g$  wherever for $r, r' \in Q$,  
$r \cap r'$ is an even sided $2$-cell.

\begin{thm}\label{knot}
Let $X$ be a cpx on $S^3$, and $G$ defined as above,
and assume $ \forall r \in Q$,  deg$(r)=2$, and $G$ has no triangles.  Then
\begin{enumerate}
\item If $G = \emptyset$ then $X$ is $4$-colorable.
\item If $G$ is an unlink then $\exists X' \approx X$ s.t. $X'$ is $5$-colorable.
\end{enumerate}
\end{thm} 
\begin{proof} 
If $|G|=0$, then $X$ ecpx on $S^3$ and so $X$ is $4$-colorable.
Let $X$ cpx on $S^3$ be given with $G \subset \cup Q$ as above and $G$ an unknot.
Let $e$ be the minimal length of all $1$-cells in $G$.
Let $M^3$ be the double cover around $G$, and $Y$ the double cover cpx corresponding to $X$.
Each cell of $X$ not meeting $G$ is produced twice in $Y$,
and the $2$-cells meeting $G$ are produced once in $Y$ with twice as many sides.
Thus $Y$ is an ecpx on $M\cong S^3$, and $Y$ may be $4$-colored.
In the projection back to $X$,  the colors of $Q$ are well defined but some of the other regions will recieve $2$-colors.
Also for $\gamma \in \pi_1(S^3 \setminus G)$, 
the cover shows $\mbox{Hol}_X(\gamma) \circ \mbox{Hol}_X(\gamma) = 1$.  
Thus if $Q$ uses colors $1$ and $2$, it can use no others as $\gamma$ is the transposition $(34)$.
So $|Q|$ is even.

Choose smooth disk $D^2$ with $D \cap \ \cup Q = \partial D$.
Choose $\epsilon < e$ so that $D_\epsilon$ is a tubular neighborhood of $D$.
Cover $D$ with fully transverse $4$-colored $3$-cells of diameter $<\epsilon/5$,
and extend this to $X'  \approx X$ with $X'^{\langle 3 \rangle}$ equal to $Q$ on $\cup Q$
and having cell diameters $< \epsilon/5$ elsewhere.
Since $\pi_1(D)=1$, we may extend the coloring on $Q$ to a local coloring on $Q \cup N_\epsilon (D)$.
$D$ is fully transverse to $X'$, and may be isotopied by moving points $<\epsilon /3$
to miss regions in color $4$.
For each region $r$ of $X'$ meeting $D$ and in color $3$,
insert $2$-cell $r\cap D$ to form cpx $X''$.
As in Theorem~\ref{make_to_ab}, $X''$ is also locally colorable around $D$,
but with the transposition in coloring $(34)$ applied to one side.
Thus all $2$-cells are even-sided and $X''$ is an ecpx on $S^3$.  
$X'$ may be colored identically to $X''$ but with color $5$ applied to the regions with inserted disks.

For $G$ an unlink, multiple disks may be chosen with $N_\epsilon(D_1, \ldots D_k)$ a
tubular neighborhood, and $X' \approx X$ with fine enough cells to avoid interference between the disks.
\end{proof}

\begin{conj}\label{unknot}
Let $X$ be a cpx on $S^3$, and $G$ defined as above,
and assume deg$(r)=2, \forall r \in Q$, and $G$ has no triangles.  Then
\begin{enumerate}
\item
 If $G$ is an unknot then $X$ is $5$-colorable.
\item If $G$ is a knot then $\exists X' \approx X$ s.t. $X'$ is $6$-colorable
\end{enumerate}
\end{conj}
The difficulty in proving (1) is that when moving $D$ to miss the color $4$ regions, one must be sure $D$ remains embedded.
For (2), further considerations are surely needed.
A Seifert surface $S$  of the knot is not simply connected so there need be no local coloring.  Furthermore
the double cover generally will not be $S^3$. It s unclear if
 Hol$_X(\gamma)$  non-trivial implies $\gamma$
intersects $S$. If all these hurdles can be overcome it appears $S$ may be isotopied
so that the double colored regions intersecting $S$ can be re-colored with colors $5$  and $6$
to separate $X'$ into a $4$-colorable ecpx.

It is interesting to consider what conditions similar to Theorem~\ref{knot} 
that bound the chromatic number of $X$
could be found when  $g$ and $G$ are more complex graphs.  
The next example refocuses on classification of $M^3$.

\begin{figure}[ht]
\begin{center}
\includegraphics*[width=1.0\textwidth,height=.3\textwidth]{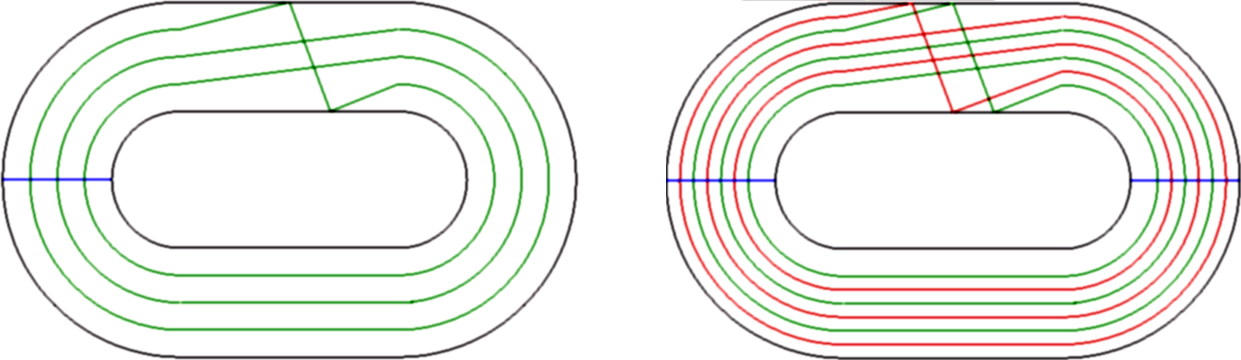}
\caption{\label{Ball_Pairs}Heegaard Diagram for Lens Space $L(3,1)$ and Decomposition into Four $B^3$'s  } 
\end{center}
\end{figure}

\begin{figure}[ht]
\begin{center}
\includegraphics*[width=.8\textwidth,height=.2\textwidth]{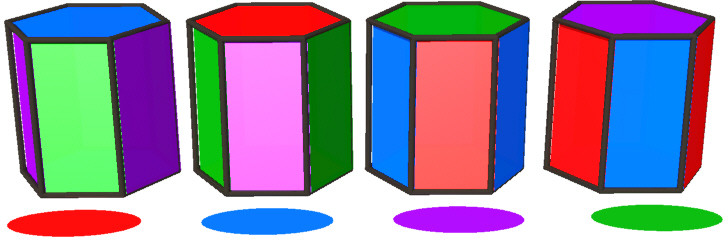}
\end{center}
\caption{\label{Four_Balls}Four $4$-Colored $B^3$ for $L(3,1)$.  Region Color shown below, External Color on Faces  } 
\end{figure}

\begin{ex}Heegaard splitting of $M^3$.
\end{ex}
Similar to the Heegaard splitting of $M^3$, Propositions~\ref{exists_cpx} gives a 
decomposition into a cpx $X$ with $4$ regions.
Let cpx $X$ on $M$ be chosen with $4$ regions $(r_1, r_2, r_3, r_4)$.
If $M^3$ is orientable, then a pairing of the regions $(r_1, r_2), (r_3, r_4)$
produces a handlebody decomposition of $M^3$.  
In Figure~\ref{Ball_Pairs}, 
the Heegaard diagram for a handlebody decomposition for the lens space $M^3 = L(3,1)$ is shown in the left panel.
On the right the solid tori have been decomposed into two $B^3$, and the resulting diagram shown.
In Figure~\ref{Four_Balls} each of the $B^3$ have been given the color shown below the regions,
and the colors of the adjoining regions shown on their faces. 
Each of the boundaries of the regions is a $3$-colored ecpx. 
 To recover the lens space the face pairing would be specified.
With this there are $3$ pairings, 
offering $3$ related Heegaard splittings.  What is the relationship between
these Heegaard splittings, and does the coloring information provide any added description?

\begin{figure}[ht]
\begin{center}
\includegraphics*[width=.4\textwidth,height=.3\textwidth]{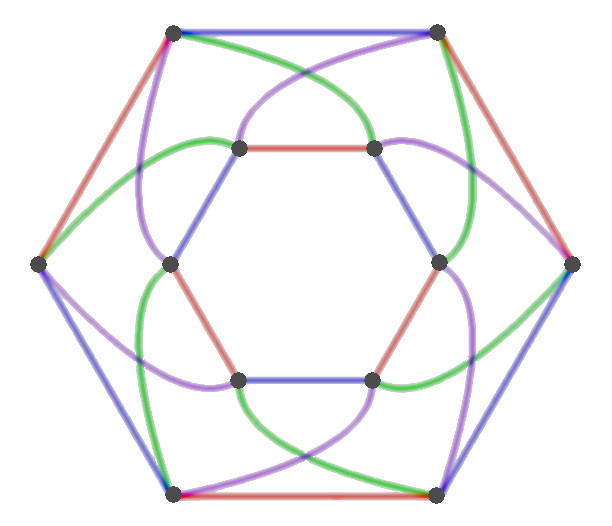}
\caption{\label{4graph}$4$-Colored Degree $4$ Graph for Lens Space $L(3,1)$   } 
\end{center}
\end{figure}

\begin{ex}Graph diagram of $M^3$.
\end{ex}
 The decomposition of $M^3$ into an ecpx $X$ with four regions gives the graph $X^{(1)}$
with degree $4$ vertices.  We may $4$-color the edges of the graph by giving each edge the remaining color distinct
from the $3$-colors of its adjacent regions.  Note for any subset of $3$-colors, the subgraph of edges
forms an ecpx on the boundary of the region of the missing $4$th color.  Thus the $3$-color 
subgraphs are all planar.  This is illustrated for the lens space $L(3,1)$ in Figure~\ref{4graph}.

Conversely suppose we are given a graph $G$ with degree $4$ vertices and $4$-colored edges.
If the $3$-color subgraphs are all planar then $G$ corresponds to a unique $M^3$.  
$G$ may be projected into $\mathbb{R}^2$ similarly to a knot diagram  with non-intersecting crossings.
How do the properties of a diagram of $G$ correspond to the properties of $M$?

The next result summarizes a key idea useful for studying $\pi_1(M)$. 
The theorem says  that homomorphisms from the fundamental group to a 
finite group are equivalent to physically represented structures on $M$.

\begin{thm}
Let $M^n$ and $j\geq 1$ be given.  Then Hol induces a bijection
\[\mbox{Hol}: \{[{\bf X}]|{\bf X} \mbox{ is a } j \mbox{-layer ecpx on } M \} \to \{f|f:\pi_1(M) \to S_{n+j} \} \]
given by $f=Hol_{\bf X}$.
\end{thm}
\begin{proof}
By Theorem~\ref{j_hol_hol}, $[{\bf X}]=[{\bf Y}] \iff \mbox{ Hol}_{\bf X}=\mbox{ Hol}_{\bf Y}$,
so Hol is well-defined on holonomy classes and injective.  By Theorem~\ref{j_make_hol_f},
$\forall f:\pi_1(M) \to S_{n+k}$,  $\exists j\mbox{-layer ecpx }{\bf X}$  on  $M$ s.t. $\mbox{Hol}_{\bf X} = f$,
so Hol is surjective.
\end{proof}

Other structures are possible.  If $M$ is orientable, $X$ ecpx on $M$,
then $X^{(1)}$ is even cyclic so we may replace alternate vertices $v$ of $X$ with a cell
$N_\epsilon(v)$ to make cpx $X'$ on $M$.  We may color the new cells with the $(n+2)$ color.
Note each edge then has $2$ distinctly colored regions at its ends and the remaining $n$ distinctly colored regions adjacent to it.
Note also that all $2$-cells are $3k$-sided, for some integer $k$.
Calling 
 any cpx $X'$ on $M$ with  $3k$-sided $2$-cells a $3$cpx,
this generalizes to give a holonomy on $S_{n+2}$ where the ends of each edge are also required to be regions of different colors.
Also Definition~\ref{Hol_Xi} is rephrased  so $Y$ and $Y'$ are like extensions.
If $X'$ is such a $3$cpx we get a map
\[\mbox{Hol}_X:\pi_1(M^n) \to A_{n+2} \]
where $A_{n+2} \subset S_{n+2}$ is the alternating group. 

For $M^n$ orientable, an $(n+1)$-colored $X$ ecpx on $M$ and an $(n+2)$-colored $X'$ $3$cpx on $M$ 
also give distinct maps $ M \to S^n$.
Let $J$ be the ecpx on $S^2$ with two vertices and three regions and $J'$ the tetrahedron.
For $3$-colored $X$ on $M^2$ there is a map $f:X \to J$
with singularity of order $k$ in the interior of each $2k$-gon region.
For $4$-colored $X'$ on $M^2$ there is a map $f':X' \to J'$
with singularity of order $k$ in the interior of each $3k$-gon region.
Other models beside $J$ and $J'$ can be used to define a holonomy.

\end{document}